\documentclass[12pt]{article} 
\usepackage[utf8]{inputenc} 
\usepackage{amsmath,amssymb}
\usepackage[margin=1in]{geometry}
\usepackage{multirow}
\usepackage[round]{natbib}
\usepackage{graphicx}
\usepackage{mathtools}
\usepackage{authblk}
\usepackage{amsthm}
\newtheorem{theorem}{Theorem}
\newtheorem{assumption}{Assumption}

\title{Asymptotic equivalence of paired Hotelling test and conditional logistic regression}
\author[1,2]{Félix Balazard}

\affil[1]{Sorbonne Universités, UPMC Univ Paris 06, CNRS}
\affil[2]{{INSERM U1169, Hôpital Bicêtre, Université Paris-Sud, felix.balazard@inserm.fr}}
\date{} 

\begin{document}
\maketitle

\section*{Abstract}
Matching, the stratification of observations, is of primary importance for the analysis of observational studies. We show that the score test of conditional logistic regression and the paired Student/Hotelling test, two tests for paired data, are asymptotically equivalent.

\textbf{Keywords:} matching, conditional logistic regression, asymptotic equivalence, paired Student test, observational studies.
\section{Introduction}

Observational studies are the main source of information in many areas of science such as epidemiology or the social and political sciences \citep{Rosenbaum2002}. Such studies aim at estimating the effect of a treatment on an outcome and they differ from randomized experiments by the absence of random allocation of treatment. This lack of randomization is due to ethical or economical considerations and as a consequence, observational studies are subject to confounding. In the typical context of an observational study, we have $n$ observations of a random outcome $Y$ and a random predictor or treatment $X$. The analysis is complicated by the possibility of confounding i.e., the influence of a third variable $Z$ on $X$ or $Y$. To take an example, suppose $Y$ is lung cancer, $X$ is yellow coloration of the teeth and $Z$ is smoker status. If we do not take into account $Z$, we will find an association between $Y$ and $X$ but if we control for $Z$, the association will disappear since both lung cancer and yellow coloration of the teeth are caused by smoking.  In a randomized experiment, randomization will ensure that the treated and non-treated group have similar values for all confounders.

Matching is an important method to control for bias in observational studies. Matching is the grouping of observations in strata. Such a grouping is necessary when the observations in a strata cannot be considered independent. For example, consider a continuous variable, such as blood pressure, which is measured before and after the beginning of a medical treatment given to a group of patients. In order to show that the treatment is effective, an incorrect approach would be to compare the mean of the group before treatment to the mean after treatment using a two-samples Student test. This naive approach relies on the unwarranted assumption of independence of the measurement of the same patient before and after the treatment. The correct approach is to compute the difference between blood pressure before and after treatment for each patient and then to test if the differences have zero mean using a Student test. This is called the paired Student test. In this simple example, the stratas are the pairs of observations of the same patient. 


In presence of confounding, matching on the values of the confounder remove bias \citep{rubin73} or when there are several confounders, other techniques based on matching can be used \citep{rosenbaum83}. Matching can also guide the data collection of an observational study. This is best illustrated by the case-control study in epidemiology. Such a study begins by gathering disease cases for example at hospitals. The recruitment of controls is less straightforward than the recruitment of cases and must try and limit confounding. One way of achieving this is to choose controls that are naturally matched to cases: for each case, controls are chosen among friends or neighbors whose --potentially unobserved-- confounders are similar to that of the case. As the case and his controls have been chosen for their similarity, they cannot be considered independent and must be grouped in a strata. 

Once a matching is obtained, adapted statistical procedures are needed and especially statistically sound testing procedures. A good reference on statistical methods for matched data is the monograph by \cite{breslow1981}. For strata limited to pairs of one case and one control (or before and after treatment) and when the predictor $X$ is a binary variable, the so-called McNemar test can be used \citep{mcnemar}. In this simple case, there are only $4$ possibilities for each pair. The data can be summarized in the following two-by-two contingency table.

\begin{table}[!ht]
\centering
\caption{\bf{McNemar test}}
\label{table1}
\begin{tabular}{llll}
                                       &                       & \multicolumn{2}{c}{$Y=1$}                          \\ \cline{3-4} 
                                       & \multicolumn{1}{l|}{} & \multicolumn{1}{l|}{$X=0$} & \multicolumn{1}{l|}{$X=1$} \\ \cline{2-4} 
\multicolumn{1}{l|}{\multirow{2}{*}{$Y=0$}} & \multicolumn{1}{l|}{$X=0$} & \multicolumn{1}{c|}{$a$} & \multicolumn{1}{c|}{$b$} \\ \cline{2-4} 
\multicolumn{1}{l|}{}                  & \multicolumn{1}{l|}{$X=1$} & \multicolumn{1}{c|}{$c$} & \multicolumn{1}{c|}{$d$} \\ \cline{2-4} 
\end{tabular}
\end{table}
In the table above, $a+d$ is the number of pairs of a case and a control concordant for the predictor $X$ and $b+c$ is the number of discordant pairs. The statistic $\xi_{mc}$ of the McNemar test depends only on the numbers of discordant pairs $b$ and $c$ and has the following form:
\begin{equation*}
\xi_{mc}=\frac{(b-c)^2}{b+c}.
\end{equation*}
Under the null hypothesis of independence between $Y$ and $X$, $\xi_{mc}$ follows the $\chi^2$ distribution asymptotically \citep{mcnemar}. This test statistic has the advantage of being computationally simple.

However, in many settings, we have access to more than one control for each case. For example, in the Isis-Diab study \citep{balazard2016}, each patient was asked to recruit two controls . A test for strata of arbitrary size and a binary predictor was proposed in \cite{MH} and is now referred to as the Cochran-Mantel-Haenszel (CMH) test. It also provides an estimate of the odd-ratio associated with the predictor. The original article contains an interesting discussion of practical aspects of the prospective study of disease.

If the predictor is not binary but continuous, appropriate tests are available if the matching is a pairing i.e., each case has exactly one control. If the vector of differences between cases and controls can be assumed to be normally distributed, a Student test can be applied on this vector to decide if its mean is zero as in our first example. The resulting test is called the 
paired Student test \citep[p180-185]{mcdonald2009}. The paired Hotelling test is the generalization of the paired Student test when there are several predictors. When the normality assumption is not respected, a non-parametric alternative is available via the Wilcoxon signed rank test.

None of the above procedures allow for analysis of continuous variables with arbitrary strata size. This was made possible by conditional logistic regression (CLR) first introduced in \cite{breslow78}. CLR is however not limited to continuous variables and can be applied in all settings of the previous tests. As CLR is based on logistic regression, it is more flexible than the previous tests and allows for multivariate analyses. 




When different statistical tests are available for the same data, it is desirable to prove  that their results are similar. For procedures based on maximum likelihood such as CLR, the three tests of significance for predictors (the Wald test, the likelihood ratio test and the score test) have been shown to be asymptotically equivalent \citep{engle1984wald}. This result of asymptotic equivalence is with respect to the null hypothesis as well as to a sequence of local alternatives.

For the analysis of matched data with a binary predictor, the available tests have been shown to be identical. When the strata are pairs, both the CMH test and the McNemar test can be applied and their statistic is identical \citep[p413]{agresti}. When the strata size is arbitrary, CMH and CLR can be applied and the CMH statistic is identical to the score test statistic of CLR \citep{Day79} and therefore if the strata are pairs, the McNemar statistic will be identical to the score test statistic of CLR. 

When the strata are pairs and the predictor (resp. predictors) is continuous, the paired Student test (resp. paired Hotelling test) and CLR can be applied. However, to the best of the author's knowledge, there are no results comparing the two procedures in that case. The objective of this note is to prove that the paired Hotelling test is asymptotically equivalent to the score test of CLR and therefore also to the Wald and likelihood ratio tests.

We start by deriving the exact form of the two test statistics in section 2. Finally, we present our result and prove it in section 3.

\section{Test statistics}
Let $Y_{i\ell}$ be the label of the $\ell$-th individual of the $i$-th stratum and $X_{i\ell}\in\mathbb{R}^p$ its continuous predictor value. The paired Hotelling test can be applied when stratas are pairs of discordant observations and therefore $\ell\in\{1,2\}$. Without loss of generality, we can assume that $Y_{i1}=1$ and $Y_{i2}=0$.

We are interested in testing the null hypothesis: \begin{equation}\mathcal{H}_0: \text{There is no association between $Y$ and $X$ }\label{null}\end{equation}
against the alternative hypothesis:
 \begin{equation*}\mathcal{H}_1: \text{There is an association between $Y$ and $X$. }\label{null}\end{equation*}

To test this null hypothesis in this setting, we have two options: the paired Hotelling test whose test statistic we denote $\xi_{\text{hot}}$ or CLR whose score test statistic we denote $\xi_{\text{sc}}$. We now derive the expressions of those two test statistics. As we will show below, the paired Hotelling test and CLR depend only on the vector of differences between pairs. We therefore denote $Z_i=X_{i1}-X_{i2}$ the difference between pairs.

\paragraph{Paired Hotelling test}
The hypotheses tested by the paired Hotelling test are expressed by the parameter  $\boldsymbol{\mu}=(\mathbb{E}[Z_j])_{j\in\{1,\hdots,p\}}$, the mean of the difference between pairs. Using this parameter, the hypothesis \eqref{null} becomes $\mathcal{H}_0 : \boldsymbol{\mu}=0$. 

Let the unbiased sample covariance be defined as:
\begin{equation*}C=\frac{1}{n-1}\sum_{i=1}^n(Z_i-\bar{Z})(Z_i-\bar{Z})^\top,\end{equation*}
where $\bar{Z}=\frac{1}{n} \sum_{i=1}^n Z_i$ and $*^\top$ is the transposition operator. The statistic of the test is: 
\begin{equation}\xi_{\text{hot}}=n\bar{Z}^\top C^{-1} \bar{Z}.\end{equation}
When $Z$ follows a centered Gaussian distribution, the distribution of $\xi_{\text{hot}}$ is called the Hotelling distribution.

\paragraph{Conditional Logistic Regression score test}
In logistic regression, the likelihood of an observation, given the parameters $\alpha\in\mathbb{R}$ and $\boldsymbol{\beta}\in\mathbb{R}^p$, equals:
\begin{equation*}\mathbb{P}(Y=y|X)=\frac{e^{y(\alpha+\boldsymbol{\beta}^\top X)}}{1+e^{\alpha+\boldsymbol{\beta}^\top X}},\end{equation*}
with $y\in\{0,1\}$.

Logistic regression can take into account stratification by having a different constant term $\alpha_i$ for each stratum. However, when the strata are small e.g., pairs, this leads to a large number of parameters and biases parameter estimation \citep[p249-251]{breslow1981}. Conditional logistic regression adresses the problem by conditioning the likelihood on the number of cases in each stratum. This eliminates the need to estimate the $\alpha_i$. In the case  of pairings, when the first observation is a case and the second is a control, the conditional likelihood of the $i$-th strata is:

\begin{align*}&\mathbb{P}(Y_{i1}=1,Y_{i2}=0|X_{i1},X_{i2},Y_{i1}+Y_{i2}=1)\\
&=\frac{\mathbb{P}(Y_{i1}=1|X_{i1}) \mathbb{P}(Y_{i2}=0|X_{i2})}{\mathbb{P}(Y_{i1}=1|X_{i1}) \mathbb{P}(Y_{i2}=0|X_{i2})+\mathbb{P}(Y_{i1}=0|X_{i1}) \mathbb{P}(Y_{i2}=1|X_{i2})}\\
&=\frac{\frac{\exp(\alpha_i+\boldsymbol{\beta}^\top X_{i1})}{1+\exp(\alpha_i+\boldsymbol{\beta}^\top X_{i1})}\times\frac{1}{1+\exp(\alpha_i+\boldsymbol{\beta}^\top X_{i2})}}{\frac{\exp(\alpha_i+\boldsymbol{\beta}^\top X_{i1})}{1+\exp(\alpha_i+\boldsymbol{\beta}^\top X_{i1})}\times\frac{1}{1+\exp(\alpha_i+\boldsymbol{\beta}^\top X_{i2})}+\frac{1}{1+\exp(\alpha_i+\boldsymbol{\beta}^\top X_{i1})}\times\frac{\exp(\alpha_i+\boldsymbol{\beta}^\top X_{i2})}{1+\exp(\alpha_i+\boldsymbol{\beta}^\top X_{i2})}}\\
\label{pair} &=\frac{\exp(\boldsymbol{\beta}^\top X_{i1})}{\exp(\boldsymbol{\beta}^\top X_{i1})+\exp(\boldsymbol{\beta}^\top X_{i2})} \qquad\text{(The $\alpha_i$ have been eliminated.)}\\
&=\frac{1}{1+\exp(\boldsymbol{\beta}^\top(X_{i2}-X_{i1}))}.\end{align*}
In the general case, when there is $k$ cases in a strata of size $m$, with the cases being the first $k$ observations, the conditional likelihood of a strata can be written:
\begin{equation*}\mathbb{P}(Y_{ij}=1\text{ for }j\leq k,Y_{ij}=0\text{ for } j>k|X_{i1},\hdots,X_{im},\sum_{j=1}^m Y_{ij}=k)=\frac{\exp(\sum_{j=1}^k \boldsymbol{\beta}^\top X_{ij})}{\sum_{J\in \mathcal{C} _{k}^{m}} \sum_{j\in J} \exp(\boldsymbol{\beta}^\top X_{ij})},\end{equation*}
where $\mathcal{C} _{k}^{m}$ is the set of all subsets of size $k$ of the set $\{1,\hdots,m\}$. 

To obtain the full conditional likelihood, we multiply over each stratum and take the $\log$ to obtain the $\log$-likelihood $L(\boldsymbol{\beta},Z)$:
\begin{equation*} L(\boldsymbol{\beta},Z)=-\sum_{i=1}^n \log(1+e^{-\boldsymbol{\beta}^\top Z_i}).\end{equation*}
The score $s(\boldsymbol{\beta},Z)$ is then defined by:
\begin{equation*}s(\boldsymbol{\beta},Z)=\frac{\partial L}{\partial \boldsymbol{\beta}}(\boldsymbol{\beta},Z)=\sum_{i=1}^n \frac{Z_i }{\exp(\boldsymbol{\beta}^\top Z_i) +1}.\end{equation*}
If $\boldsymbol{\beta}=\boldsymbol{\beta}_0$, the covariance matrix of $s(\boldsymbol{\beta}_0,Z)$ is Fisher's information matrix:
\begin{equation*}\mathcal{I}(\boldsymbol{\beta}_0)=\left(\mathbb{E}[-\frac{\partial^2 L}{\partial \boldsymbol{\beta}_j\partial\boldsymbol{\beta}_k}(\boldsymbol{\beta}_0,z)|\boldsymbol{\beta}_0 ]\right)_{j,k\in\{1,\hdots,p\}}=\left(\sum_{i=1}^n \frac{Z_{ij}Z_{ik} \exp(\boldsymbol{\beta}_0^\top Z_i)}{(\exp(\boldsymbol{\beta}_0^\top Z_i)+1)^2}\right)_{j,k\in\{1,\hdots,p\}}.\end{equation*}
To test the hypothesis $\mathcal{H}_0: \boldsymbol{\beta}=\boldsymbol{\beta}_0$, the score test statistic is:
\begin{equation*}\xi_{\text{sc}}=s(\boldsymbol{\beta}_0,Z)^\top\mathcal{I}(\boldsymbol{\beta}_0,z)^{-1}s(\boldsymbol{\beta}_0,Z)\end{equation*}

Expressed using $\boldsymbol{\beta}$, hypothesis \eqref{null} becomes $\mathcal{H}_0: \boldsymbol{\beta}=0$. As shown above, the null hypothesis of the paired Hotelling test is expressed using $\boldsymbol{\mu}$, the mean of the difference between pairs, rather than $\boldsymbol{\beta}$, the regression coefficient. There is no general correspondence between the two parameters. However, the two null hypotheses that are of interest in practice are the same: $\boldsymbol{\mu}=0 \iff \boldsymbol{\beta}=0$.

For $\boldsymbol{\beta}=0$, the score, Fisher's information matrix and the test statistic become:
\begin{equation*}s(0,Z)=\frac{n}{2}\bar{Z},\qquad \mathcal{I}(0,Z)=\frac{1}{4}\sum_{i=1}^n Z_i Z_i^\top\end{equation*}
and
\begin{equation}\xi_{\text{sc}}=n\bar{Z}^\top\left(\frac{1}{n}\sum_{i=1}^n Z_iZ_i^\top\right)^{-1}\bar{Z}.\end{equation}

Now that we have expressed analytically $\xi_{\text{hot}}$ and $\xi_{\text{sc}}$, we state our result in the following section.

\section{Asymptotic equivalence}

In this section, we introduce the sequence of local alternatives and state the result of asymptotic equivalence between the two tests.

We adopt the framework of  \cite{engle1984wald} for asymptotic equivalence. The motivation for considering a sequence of local alternatives is that any reasonable test will have the right size and will reject a fixed alternative when the number of observations becomes large. The sequence of local alternatives considers deviations from the null that approach the null as sample size increases. This allows to compare tests more precisely.

We model the sequence of local alternatives by a triangular array of observations
$(Z_i^{(n)})_{n\in\mathbb{N},i\in\{1,\hdots,n\}}$ that respects the following assumption. Let $\boldsymbol{\delta}\in\mathbb{R}^p$ and $\Sigma$ be a positive definite matrix of dimension $p$.
\begin{assumption} \label{triangular}
For all $n\in \mathbb{N}$, $i\in\{1,\hdots,n\}$, we have $\mathbb{E}[Z_i^{(n)}]=\frac{\boldsymbol{\delta}}{\sqrt{n}}$ and $\text{\textnormal{Cov}}[Z_i^{(n)}]=\Sigma$. In addition, $W_i^{(n)}=Z_i^{(n)}-\frac{\boldsymbol{\delta}}{\sqrt{n}}$ are independent and identically distributed.\end{assumption}
The triangular array of \textbf{Assumption \ref{triangular}} is the sequence of local alternatives. The parameter $\boldsymbol{\delta}$ controls the deviation from the null. The null hypothesis is the special case of the sequence of local alternatives when $\boldsymbol{\delta}=0$. The $\sqrt{n}$ in the denominator of the deviation shrinks the deviation towards the null as the sample size increases. This particular power of $n$ ensures that the test statistics are well-behaved as they converge in distribution.

We denote $\xi_{\text{sc},n}$ (resp. $\xi_{\text{hot},n}$) the statistic of the score test (resp. of the paired Hotelling test) associated with $Z^{(n)}$, the $n$-th line of the triangular array of \textbf{Assumption \ref{triangular}}. In the same fashion as for the Student distribution and the $\chi^2$ distribution, the Hotelling distribution with degrees of freedom $p$ and $n-1$ and the $\chi^2_p$ distribution are different for finite samples but they are the same asymptotically. Therefore,  $\xi_{\text{sc},n}-\xi_{\text{hot},n} \xrightarrow{P} 0 $ is enough to guarantee asymptotic equivalence. Our theorem is more precise as it includes the convergence rate and limiting distribution of this difference. We note the convergence in distribution $d$.

\begin{theorem}Under the null or a sequence of local alternatives, the paired Hotelling test and the score test of CLR are asymptotically equivalent. More precisely, under \textbf{Assumption \ref{triangular}}, we have: $$n(\xi_{\text{\textnormal{sc}},n}-\xi_{\text{\textnormal{hot}},n})\xrightarrow{d} K$$ where  $K=(\boldsymbol{\delta}+V)^\top\Sigma^{-1}\left(\Sigma-(\boldsymbol{\delta}+V)(\boldsymbol{\delta}+V)^\top\right)\Sigma^{-1}(\boldsymbol{\delta}+V)^\top$ with $V\thicksim\mathcal{N}(0,\Sigma)$ and $\mathcal{N}$ refers to the normal distribution. \end{theorem}

\begin{figure}\begin{center}\includegraphics[scale=0.6]{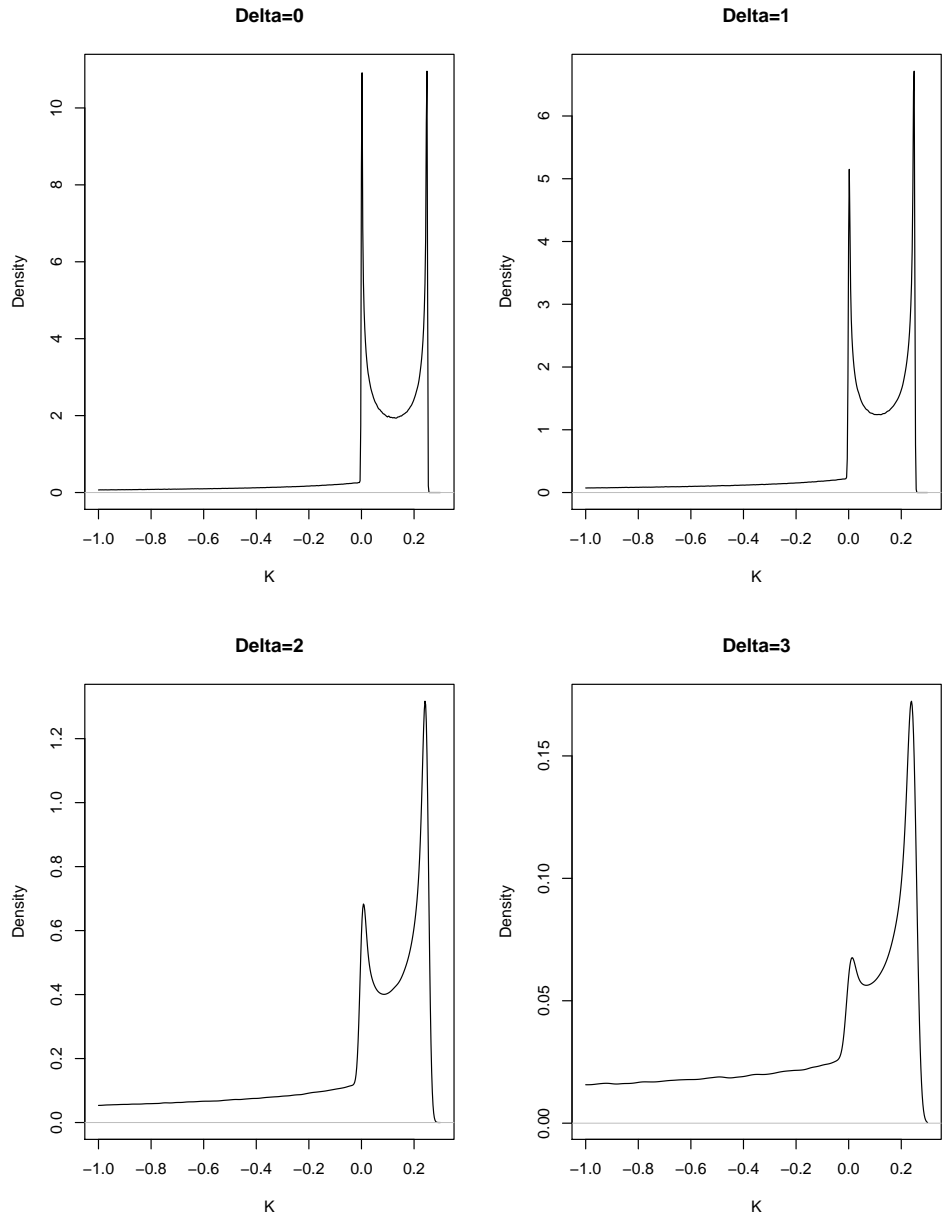}\caption{\textbf{Distribution of K.} Distribution of the limit variable for $p=1$, $\Sigma=\sigma^2=1$ and $\boldsymbol{\delta}\in\{0,1,2,3\}$. }\label{distributionK} \end{center}\end{figure}

Examples of distribution of $K$ for $p=1$, $\Sigma=\sigma^2=1$ and different values of $\boldsymbol{\delta}$ are shown in figure \ref{distributionK}. We see that for small values of $\boldsymbol{\delta}$, K is concentrated between $0$ and $1/4$. As $\boldsymbol{\delta}$ grows, the distribution shifts to large, negative values.

For the practitioner, the result of \textbf{Theorem \ref{triangular}} will not affect the choice of test but ensures that similar conclusions are reached regardless of that choice. The result also implies that CLR can be considered a generalization of the paired Hotelling or Student test when the conditions for the latter do not apply.

\paragraph{Proof of Theorem 1}
To facilitate the notation, we omit the dependence in $n$ in most quantities, writing $Z$ instead of $Z^{(n)}$. The quantity we want to estimate is:
\begin{equation}n(\xi_{\text{sc},n}-\xi_{\text{hot},n})=n^2\bar{Z}^\top(\tilde{\mathcal{I}}^{-1}-C^{-1})\bar{Z}=n^2\bar{Z}^\top\tilde{\mathcal{I}}^{-1}(C-\tilde{\mathcal{I}})C^{-1}\bar{Z},\label{diff}\end{equation}
where $\tilde{\mathcal{I}}=\frac{1}{n}\sum_{i=1}^n Z_i Z_i^\top$.

The Lindeberg-Feller Central Limit Theorem for triangular arrays \citep[Theorem 3.2 p58]{hall_heyde} applied on the $W_i^{(n)}$ gives:
\begin{equation}\label{CLT} \sqrt{n}\bar{Z}\xrightarrow{d} (\boldsymbol{\delta} +V),\end{equation}
with $V\thicksim\mathcal{N}(0,\Sigma)$.
The weak law of large numbers for triangular arrays implies:
\begin{equation}C\xrightarrow{P} \Sigma \label{LLN} .\end{equation}
Combining \eqref{CLT} and \eqref{LLN}, we see:
\begin{equation}\label{LLN2} \tilde{\mathcal{I}}=\frac{1}{n} \sum_{i=1}^n Z_i Z_i^\top=\frac{1}{n} \sum_{i=1}^n (Z_i-\bar{Z})(Z_i-\bar{Z})^\top+\bar{Z}\bar{Z}^\top \xrightarrow{P} \Sigma\end{equation}
and, using the second identity in \eqref{LLN2}, we obtain:
\begin{equation*}n(C-\tilde{\mathcal{I}})=C-n\bar{Z}\bar{Z}^\top\xrightarrow{d}\Sigma-(\boldsymbol{\delta}+V)(\boldsymbol{\delta}+V)^\top.\end{equation*}
Applying  Slutsky's lemma to all of the above, we conclude:
\begin{equation*}n(\xi_{\text{sc},n}-\xi_{\text{hot},n})\xrightarrow{d}K\end{equation*}
where $K=(\boldsymbol{\delta}+V)^\top\Sigma^{-1}\left(\Sigma-(\boldsymbol{\delta}+V)(\boldsymbol{\delta}+V)^\top\right)\Sigma^{-1}(\boldsymbol{\delta}+V)^\top$ with $V\thicksim\mathcal{N}(0,\Sigma)$.
This concludes the proof of \textbf{Theorem 1}.
\begin{flushright}$\qed$\end{flushright}

\section*{Acknowledgements}
The author thanks Gérard Biau for his supervision and his many comments on the manuscript.
Funding: The author acknowledges a PhD grant from École Normale Supérieure.

\section*{References}

\bibliography{bib}
\end{document}